\documentclass[12pt, reqno]{amsart}
\usepackage{amsmath, amstext, amsbsy, amssymb, amscd}

\newtheorem{theorem}{Theorem}[section]
\newtheorem{lemma}[theorem]{Lemma}
\newtheorem{proposition}[theorem]{Proposition}
\newtheorem{corollary}[theorem]{Corollary}
\theoremstyle{definition}

\theoremstyle{remark}
\newtheorem{remark}[theorem]{Remark}

\numberwithin{equation}{section}

{\vskip-\lastskip\medskip
  \noindent
  {\em #1.}\enspace
  }%
{\qed\par\medskip
  }

\begin{document}
\title[ On
       threefolds without  nonconstant regular functions ] 
              { On
       threefolds without  nonconstant regular functions } 

\author[Jing  Zhang]{Jing  Zhang}
\address{Department of Mathematics, University of 
Missouri, Columbia, MO
65211, USA}
\email{zhangj@math.missouri.edu}


\begin{abstract}

We consider smooth threefolds $Y$ defined over $\Bbb{C}$ with  
$H^i(Y, \Omega^j_Y)=0$
for all $j\geq 0$, $i>0$. 
Let $X$ be a smooth projective threefold
containing $Y$
and $D$ be the boundary divisor with support $X-Y$.
  We are interested in the following question:
What  geometry  information  of  $X$  can be  obtained 
from the regular function  information on $Y$? 
 Suppose that the boundary $X-Y$ is a smooth projective 
surface. 
In this paper, we analyse two   different cases,
i.e., 
there are no nonconstant regular functions on $Y$ or there are
lots of regular functions on $Y$.  More precisely,     
if  $H^0(Y,  {\mathcal{O}}_Y)=\Bbb{C}$,
 we prove that  
$\frac{1}{2}(c_1^2+c_2)\cdot D=\chi({\mathcal{O}}_D)\geq 0$.
In particular, if the line bundle  ${\mathcal{O}}_D(D)$ is not torsion,
then 
$q=h^1(X, {\mathcal{O}}_X)=0$,
$\frac{1}{2}(c_1^2+c_2)\cdot D=\chi({\mathcal{O}}_D)=0$, 
 $\chi({\mathcal{O}}_X) >0$  and   $K_X$ is not nef. 
If  there is  a positive  constant $c$ such that 
$h^0(X, {\mathcal{O}}_X(nD))\geq c n^3$ for all sufficiently large $n$
(we say that $D$ is big or the $D$-dimension of $X$ is 3) and $D$ has no 
exceptional curves,  then $|nD|$ is base point free for $n\gg 0$. 
Therefore $Y$ is affine if $D$ is big.  
\end{abstract}

\maketitle
\date{}
\section{Introduction}

  Let $Y$ be an algebraic manifold (i.e., an irreducible 
smooth algebraic variety defined over $\Bbb{C}$)
with $H^i(Y, \Omega^j_Y)=0$   for all $j\geq 0$ and $i>0$,
 where $\Omega^j_Y$  is the sheaf of regular  $j$-forms.
We want to understand what  $Y$  is. This question was raised by
J.-P. Serre for complex manifolds \cite{Se}. 
Since $Y$ is not compact,  for 
any analytic or algebraic coherent sheaf ${\mathcal{F}}$ on $Y$,
we have  $H^3(Y, {\mathcal{F}})=0$  \cite{Siu1, Siu2,  Zh1}. 
 So $Y$ contains
no complete surfaces \cite{NS,  Zh1}. If $Y$ has non-constant regular
functions, we know that it contains no complete curves \cite {Zh1}.
Let $X$ be a smooth completion of $Y$, then the complement of $Y$
in $X$ is connected \cite{Zh1}. Suppose that the boundary $X-Y$ 
is of pure codimension 1 and is the support   of an effective divisor 
$D$ with simple normal crossings.  We consider the 
$D$-dimension of $X$ in order 
to understand $Y$. 
The notion of  $D$-dimension is due to Iitaka
 (\cite{I1},Lecture 3 or  \cite{Uen1}, Chapter 2). 
It measures that how many  regular  functions there are on $Y$. 
If for all integers  $m> 0$ we have
     $H^0(X, {\mathcal{O}}_X(mD))=0$, then we define 
     the  $D$-dimension of $X$, denoted by $\kappa (D, X)$, to be $-\infty$.
     If  $h^0(X, {\mathcal{O}}_X(mD))\geq 1$ for some $m$, 
     choose a basis $\{f_0, f_1, \cdot \cdot\cdot, f_n\}$
     of the linear space 
     $H^0(X, {\mathcal{O}}_X(mD))$, it defines a rational 
     map 
     $\Phi _{mD}$
     from $X$ to the projective space 
     ${\Bbb{P}}^n$ by sending a point $x$ on $X$ to
     $(f_0(x), f_1(x), \cdot \cdot\cdot, f_n(x))$ in ${\Bbb{P}}^n$.  
     Then we define
     $\kappa (D, X)$ to be the maximal dimension of the images
      of the rational map  $\Phi _{mD}$, i.e., 
      $$ \kappa (D, X)= \max_m\{\dim (\Phi _{mD}(X))\}. 
      $$
Let $K_X$ be the canonical divisor of $X$, then the Kodaira 
    dimension of $X$ is  the $K_X$-dimension of $X$, denoted by 
    $\kappa(X)$, i.e.,  
    $$\kappa(X)=\kappa(K_X, X). 
    $$  
When $\kappa(D, X)=$dim$X$, we say that $D$ is big.

In our case, since $D$ is effective, the $D$-dimension
 $\kappa(D, X)\neq -\infty$. 
In \cite{Zh1, Zh2},
we prove that $\kappa(D, X)$  can be 1 and in this case, 
we have a surjective morphism from $Y$ to a smooth affine curve 
$C$ such that every fibre $S$   satisfies the same vanishing 
condition, i.e., $H^i(S, \Omega^j_S)=0$
for all $j\geq 0$ and $i>0$.  
We obtained that
 the Kodaira dimension 
of $X$  is  $-\infty$  and $q(X)=h^1(X, {\mathcal{O}}_X)$
can be any non-negative integer. In particular, all
smooth fibres are of the same type, 
i.e., type (2) or type (3) open surface   in the following
Theorem 2.6 \cite{Ku}.   
For the existence of the non-affine and non-product 
threefolds  with  vanishing   Hodge   cohomology , see \cite{Zh2}.

 In \cite{Zh3}, we proved that  
$\kappa(D, X)\neq 2$ and if $\kappa(D, X)=3$, then $Y$
is birational to Spec$\Gamma (Y, {\mathcal{O}}_Y )$. 
Furthermore, if $Y$ is regularly separable, i.e., 
if 
for any two distinct points $y_1$ and $y_2$ on $Y$, there is a regular
function $f$ on $Y$ such that $f(y_1)\neq f(y_2)$, 
 then 
$Y$ is affine. Regular separability implies 
$\kappa(D, X)=3$. 
We want to know whether $Y$ is affine 
if  $H^i(Y, \Omega^j_Y)=0$ and  $\kappa(D, X)=3$. 
If we can prove the converse, that is, 
$\kappa(D, X)=3$ implies regular separability, 
then  $Y$ is affine.  
It is sufficient to prove that if
$H^i(Y, \Omega^j_Y)=0$
for all $j\geq 0$ and $i>0$ and $\kappa(D, X)=3$,
then any non-constant regular function on $Y$ defines
an affine surface.
Another possible approach is to prove that $|nD|$ is base point free 
for some $n>0$. Since $Y$ contains no complete 
curves, $Y$ is affine (\cite{H2}, Chapter 2, Proposition 2.2). 
In this paper, we will prove  that  when 
$D$ is smooth, irreducible and contains no exceptional curves,
then $Y$ is affine.

\begin{theorem}Let $Y$ be an irreducible  smooth threefold  contained 
in a smooth projective
threefold $X$ such that  $H^i(Y, \Omega^j_Y)=0$
for all $j\geq 0$ and $i>0$. Suppose that the complement  $D=X-Y$ is a smooth 
projective 
surface without exceptional curves and $\kappa(D, X)=3$,  then 
$|nD|$ is base point free for all $n\gg 0$.
\end{theorem}

\begin{corollary} Under the condition of Theorem 1.1, 
$Y$ is affine
if and only if $\kappa(D, X)=3$. 
\end{corollary}

The most mysterious case is $\kappa(D, X)=0$. 
It is hard to understand  because 
 we cannot construct a fibre space by the divisor $D$. 
However, in order to keep track of $Y$ and its cohomology, 
we have to use this boundary divisor
to define the map. So we cannot apply 
Iitaka's  fibration  or Mori's construction.
The situation on any normal and complete  surface
is much better. Any effective divisor on a normal complete 
surface has a unique Zariski decomposition 
and any two divisors have an intersection number \cite{Sa2}.
We even can get satisfied information of   the surface
if we know the numerical type of the  divisor 
\cite{Sa1}. When dim$Y=3$ and the $D$-dimension is
0 or 3, we do not have  Zariski decomposition \cite{C} and   
a  good method to check whether a divisor is nef.       
When $\kappa(D, X)=0$ and $D$ is smooth and 
irreducible, we can reduce the problem to surface case.

\begin{theorem} Let $Y$ be a smooth threefold contained in a 
smooth projective
threefold $X$ such that   $H^i(Y, \Omega^j_Y)=0$
for all $j\geq 0$ and $i>0$ and the complement  $D=X-Y$ is a smooth 
projective 
surface.   If  $Y$
has no nonconstant regular functions, then 
$\frac{1}{2}(c_1^2+c_2)\cdot D=\chi({\mathcal{O}}_D)\geq 0$.
In particular, if the line bundle  ${\mathcal{O}}_D(D)$ is not torsion,
then  $q=h^1(X, {\mathcal{O}}_X)=0$,
$\frac{1}{2}(c_1^2+c_2)\cdot D=\chi({\mathcal{O}}_D)=0$,
  $\chi({\mathcal{O}}_X) >0$
and $K_X$ is not nef.  
\end{theorem}

We  organize the paper as follows. 
In Section 2, we will present some preparation.
 We will prove the two 
theorems in Section 3. 

Let us mention some open problems for the Steinness of $Y$. 
When dimension of $Y$ is 2, the type (3) open surfaces in Theorem 2.6
is a mystery. When dimension is 3 and $\kappa(D, X)=1$, we do not know
whether the nonaffine, nonproduct example in \cite{Zh2} is  Stein. 
Hartshorne asked the following question (\cite{H2}, page 235):
Remove an irreducible curve from a smooth complete surface, what is 
the condition for the open surface to be Stein? See \cite{N, Ued}
for recent progress. Unfortunately, their approach cannot be applied to
type (3) surface   because  the boundary has 9 components. 
We can ask  three dimensional analogue of Hartshorne's question:
What is the necessary and sufficient condition 
of a smooth  threefold to be Stein but not affine? 
Are the two conditions  $H^i(Y, \Omega^j_Y)=0$
for all $j\geq 0$ and $i>0$ and  $\kappa(D, X)\leq 1$ sufficient? 
The known result is that  
there exist  algebraic Stein  threefolds  $Y$  with  $H^i(Y, \Omega^j_Y)=0$
for all $j\geq 0$,  $i>0$ and  $\kappa(D, X)=1$ \cite{Zh1}.
However, if  $\kappa(D, X)=0$, 
we do not know whether 
there exists an algebraic Stein threefold
with   $H^i(Y, \Omega^j_Y)=0$  whereas such surface exists 
(\cite{H2}, Chapter VI, Section 3 or \cite{Ku}).

\noindent 
{\bf{Acknowledgments}}  \quad I would like to thank 
the following 
      professors for  helpful discussions: Steven Dale Cutkosky,
Dan  Edidin, 
N.Mohan Kumar, Zhenbo Qin  and Qi Zhang.

\section{Preparation}

\begin{lemma}{\bf[Goodman, Hartshorne]}  Let $V$ be a scheme and  
$D$ be an effective Cartier divisor on  $V$. Let $U=V-$Supp$D$ and $F$ 
be any coherent sheaf on $V$, 
then for every $i\geq 0,$ 
$$\lim_{{\stackrel{\to}{n}}}
H^i(V, F\otimes {\mathcal{O}}(nD)) \cong  H^i(U,  F|_U).
$$
\end{lemma}

\begin{lemma} Let $I=\{i\}$ be a direct system of indices. 
Let $\{{\mathcal{F}}_i, f_{ji}\}$, $\{{\mathcal{G}}_i, g_{ji}\}$
and $\{{\mathcal{H}}_i, h_{ji}\}$  be direct  system  of 
coherent sheaves  indexed by $I$ over a topological space $X$.
If for all $i\in I$, there are short exact sequences 
$$  0
\longrightarrow {\mathcal{F}}_i
\longrightarrow  {\mathcal{G}}_i
\longrightarrow  {\mathcal{H}}_i
\longrightarrow  0
$$
and the commutative diagram  for  all 
$j\geq i$
\[
  \begin{array}{ccccccccc}
0\longrightarrow{\mathcal{F}}_i&
\longrightarrow   {\mathcal{G}}_i&
\longrightarrow  {\mathcal{H}}_i  &\longrightarrow 0\\
\quad \quad   \Big\downarrow\vcenter{%
        \rlap{$\scriptstyle{f_{ji}}$}} &
\quad\Big\downarrow\vcenter{%
        \rlap{$\scriptstyle{g_{ji}}$}} & 
\quad \Big\downarrow\vcenter{%
        \rlap{$\scriptstyle{h_{ji}}$}}  & \\
0\longrightarrow{\mathcal{F}}_j &\longrightarrow {\mathcal{G}}_j & 
 \longrightarrow {\mathcal{H}}_j &
\longrightarrow 0, 
\end{array}
\]
then we have exact sequence
$$ 0
\longrightarrow
 \lim_{{\stackrel{\to}{i}}}   {\mathcal{F}}_i
   \longrightarrow     \lim_{{\stackrel{\to}{i}}}  {\mathcal{G}}_i
\longrightarrow 
        \lim_{{\stackrel{\to}{i}}}  
 {\mathcal{H}}_i
\longrightarrow  0.
$$
\end{lemma}
{\it Proof}. By the assumption, for any point $x\in X$, we have 
short exact sequence on stalks 
$$  0
\longrightarrow ({\mathcal{F}}_i)_x
\longrightarrow  ({\mathcal{G}}_i)_x
\longrightarrow  ({\mathcal{H}}_i)_x
\longrightarrow  0
$$
and the commutative diagram  for  all 
$j\geq i$
\[
  \begin{array}{ccccccccc}
0\longrightarrow ({\mathcal{F}}_i)_x&
\longrightarrow   ({\mathcal{G}}_i)_x&
\longrightarrow  ({\mathcal{H}}_i)_x  &\longrightarrow 0\\
\quad \quad   \Big\downarrow\vcenter{%
        \rlap{$\scriptstyle{({f_{ji}})_x}$}} &
\quad\Big\downarrow\vcenter{%
        \rlap{$\scriptstyle{{(g_{ji}}})_x$}} & 
\quad \Big\downarrow\vcenter{%
        \rlap{$\scriptstyle{{(h_{ji}})}_x$}}  & \\
0\longrightarrow ({{\mathcal{F}}_j})_x &\longrightarrow ({{\mathcal{G}}_j})_x & 
 \longrightarrow ({\mathcal{H}}_j)_x &
\longrightarrow 0. 
\end{array}
\]
By Ueno, Algebraic Geometry 2 (\cite{Uen2}, Page 10),
 we have the following 
exact sequences  of abelian groups on stalks
$$ 0
\longrightarrow
 \lim_{{\stackrel{\to}{i}}}   ({{\mathcal{F}}_i})_x
   \longrightarrow     \lim_{{\stackrel{\to}{i}}}  ({{\mathcal{G}}_i})_x
\longrightarrow 
        \lim_{{\stackrel{\to}{i}}}  
 ({{\mathcal{H}}_i})_x
\longrightarrow  0.
$$
Since  direct limits commute with each other (\cite{Br}, page 20),
we have 
$$ \lim_{{\stackrel{\to}{i}}}   ({{\mathcal{F}}_i})_x
=\lim_{{\stackrel{\to}{i}}}
\lim_{{\stackrel{\to}{x\in U}}}{{\mathcal{F}}_i}
=\lim_{{\stackrel{\to}{x\in U}}}
\lim_{{\stackrel{\to}{i}}}   {{\mathcal{F}}_i}
=(\lim_{{\stackrel{\to}{i}}}   {{\mathcal{F}}_i})_x.
$$
The Lemma follows from 
$$ 0
\longrightarrow
( \lim_{{\stackrel{\to}{i}}}   {\mathcal{F}}_i)_x
   \longrightarrow    ( \lim_{{\stackrel{\to}{i}}}  {\mathcal{G}}_i)_x
\longrightarrow 
       ( \lim_{{\stackrel{\to}{i}}}  
 {\mathcal{H}}_i)_x
\longrightarrow  0.
$$ 
\begin{flushright}
 Q.E.D. 
\end{flushright}

\begin{lemma}  Let  $Y$ be a smooth variety contained in a smooth
projective variety $X$ such that the complement $X-Y$ is compact and of 
pure codimension 1.  Let  $D$  be  any effective  divisor  with 
support  $X-Y$. Then we  have the following two exact sequences
$$ 0
\longrightarrow
 \lim_{{\stackrel{\to}{n}}}   {\mathcal{O}}_X(nD)
   \longrightarrow     \lim_{{\stackrel{\to}{n}}}  {\mathcal{O}}_X((n+1)D)
\longrightarrow 
        \lim_{{\stackrel{\to}{n}}}  
 {\mathcal{O}}_D((n+1)D)
\longrightarrow  0
$$
and for all $j>0$
$$0
\longrightarrow
 \lim_{{\stackrel{\to}{n}}} \Omega^j_X(nD)
\longrightarrow     \lim_{{\stackrel{\to}{n}}} 
\Omega^j_X((n+1)D)
\longrightarrow     \lim_{{\stackrel{\to}{n}}} 
\Omega^j_X((n+1)D)|_D
\longrightarrow  0.
$$
\end{lemma}
{\it Proof}.  For any positive integers $n$ and $m$, we have the following 
commutative diagram 
\[
  \begin{array}{ccccccccc}
0\longrightarrow{\mathcal{O}}_X(nD)&
{\stackrel{f}{\longrightarrow}} {\mathcal{O}}_X((n+1)D)&
{\stackrel{r}{\longrightarrow}} {\mathcal{O}}_D((n+1)D)  &\longrightarrow 0\\
\quad \quad   \Big\downarrow\vcenter{%
        \rlap{$\scriptstyle{i}$}} &
\quad\Big\downarrow\vcenter{%
        \rlap{$\scriptstyle{i}$}} & 
\quad \Big\downarrow\vcenter{%
        \rlap{$\scriptstyle{h}$}}  & \\
0\longrightarrow{\mathcal{O}}_X((n+m)D) &
{\stackrel{f}{\longrightarrow}} 
{\mathcal{O}}_X((n+m+1)D) & 
 {\stackrel{r}{\longrightarrow}}  {\mathcal{O}}_D((n+m+1)D) &
\longrightarrow 0, 
\end{array}
\]
where the first map $f$ in each arrow is defined by the local 
defining function of $D$, $r$ 
is the restriction map, 
 the vertical map $i$ is the natural
embedding map "1" and the last vertical map $h$ is defined as follows.
For any nonzero element $s$ in $O_D((n+1)D)$ (locally), there is an $t$
in $ O_X((n+1)D)$ (locally) such that $r(t)=s$. Then we define $h(s)$ to be
$r(i(t))=r(t)$. Since  $i\neq f$, if $t$ is not zero on $D$,
then $i(t)$ does not sit in the image of $O_X((n+m)D)$ under the map $f$. 
Notice that the vertical  map $i$ defines the direct limit. 
Then the 
first exact sequence 
 is an immediate consequence of Lemma 2.2. 
The second exact sequence can be similarly proved. 
\begin{flushright}
 Q.E.D. 
\end{flushright}

\begin{lemma} Let $Y$   be  a smooth threefold with 
 $H^i(Y, \Omega^j_Y)=0$
for all $j\geq 0$ and $i>0$.  Let  $X$   be   a smooth completion
of $Y$  such that the complement $X-Y$ is compact and  of pure codimension 1.
Let $D$ be any effective divisor on $X$ with support $X-Y$, then
for all $j\geq 0$ and $i>0$,
we have 
$$ \lim_{{\stackrel{\to}{n}}} H^i(D,  \Omega^j_X(nD)|_D)=0.  
$$
\end{lemma}
{\it Proof}. Since  $H^i(Y, \Omega^j_Y)=0$
for all $j\geq 0$ and $i>0$, by Lemma 2.1, we have 
$$\lim_{{\stackrel{\to}{n}}} H^i(X,  \Omega^j_X(nD))=0.
$$
Since  the direct limit commutes with cohomology
(\cite{H1}, Chapter III, Proposition 2.9), the Lemma follows
from the short exact sequences  in Lemma 2.3.    
\begin{flushright}
 Q.E.D. 
\end{flushright}

The following lemma is not new. We already proved it 
in our previous paper \cite{Zh1}. Since we frequently use
the argument in the proof of the theorems, We include a proof here
for completeness. 
\begin{lemma} Under the condition of Lemma 2.4, 
$H^3(X, {\mathcal{O}}_X(nD))=0$ for sufficiently large 
$n$.
\end{lemma}
{\it Proof}.  From the exact sequence 
$$ 0
\longrightarrow
  {\mathcal{O}}_X(nD)
   \longrightarrow      {\mathcal{O}}_X((n+1)D)
\longrightarrow 
 {\mathcal{O}}_D((n+1)D)
\longrightarrow  0,
$$
we have surjective map from $H^3(X,{\mathcal{O}}_X(nD))$
to $H^3(X,{\mathcal{O}}_X((n+1)D))$ for all $n\geq 0$. 
By Lemma 2.1, we have 
$$  \lim_{{\stackrel{\to}{n}}} H^3(X,{\mathcal{O}}_X(nD))=0.
$$
These two conditions imply the vanishing of the third cohomology 
for large n. 
\begin{flushright}
 Q.E.D. 
\end{flushright}

\begin{theorem}{\bf [Mohan Kumar]}   Let
$Y$ be a smooth algebraic surface over $\Bbb{C}$
with   $H^i(Y, \Omega^j_Y)=0$   for all $j\geq 0$ and $i>0$,
then  $Y$ is one of the following 
 
      (1) $Y$ is affine.

      (2) Let $C$ be an elliptic curve and $E$ the unique nonsplit 
      extension of $\mathcal{O}$$_C$ by itself.  
      Let ${X=\Bbb{P}}_C(E)$ and  $D$ be the canonical section, then $Y=X-D$.

       (3) Let $X$ be a projective rational surface with an effective 
       divisor $D=-K$ with $D^2=0$, $\mathcal{O}$$(D)|_D$ be nontorsion and 
       the dual graph of $D$ be $\tilde{D}_8$ or $\tilde{E}_8$, then $Y=X-D$.
\end{theorem}

In the above theorem, when $Y$ is affine, 
we can choose $D$ such  that  $D$ is ample 
(\cite{H2}, Theorem 4.2, page 69), so
 $\kappa(D, X)=2$. 
If $Y$ is not  affine, then  $\kappa(D, X)=0$ by Lemma 1.8 \cite{Ku}. 
\begin{theorem}{\bf [Cutkosky, Srinivas]} Let $k$ be a field
 of characteristic  0. Let $X$ be a normal surface, proper over $k$,
and let $D$ be an effective Cartier divisor on $X$. Then, for sufficiently
large $n$,  
$$  h^0(X, {\mathcal{O}}_X(nD))=P(n)+\lambda(n), 
$$
where $P(n)$ is a quadratic polynomial and the function 
$\lambda(n)$ is periodic.
\end{theorem}

For the proof of the following Iitaka's theorem, see Lecture 3 \cite{I1}
or Theorem 8.1 \cite{Uen1}.  

\begin{theorem}{\bf[Iitaka]} 
Let $X$ be a normal projective variety and let 
$D$ be an effective divisor on $X$. There exist two positive
numbers $\alpha$ and $\beta$ such that for all sufficiently 
large $n$ we have
$$  \alpha n^{\kappa(D, X)}
\leq h^0(X, {\mathcal{O}}_X(nD))
\leq \beta n^{\kappa(D, X)}. 
$$

\end{theorem}

\section{Proof of the Theorems}

From now on, we will fix the notations as follows. $Y$   is
 a smooth threefold with 
 $H^i(Y, \Omega^j_Y)=0$
for all $j\geq 0$ and $i>0$.  $X$ is  a smooth completion
of $Y$  such that the complement $X-Y=D$ is a smooth
 projective surface.

{\bf Proof of Theorem 1.1.} Since  $\kappa(D, X)=3$ and $D$ is effective,
by Theorem 2.8,  there is a constant $c>0$  such that for sufficiently large $n$,
$h^0(X, {\mathcal{O}}_X(nD))> c n^3$.  
From the short exact sequence
$$0\longrightarrow 
{\mathcal{O}}_X(nD)
\longrightarrow
{\mathcal{O}}_X((n+1)D) 
\longrightarrow
{\mathcal{O}}_D((n+1)D) 
\longrightarrow
0, 
$$
we have 
$$0\longrightarrow 
H^0({\mathcal{O}}_X(nD))
\longrightarrow
H^0({\mathcal{O}}_X((n+1)D)) 
\longrightarrow
H^0({\mathcal{O}}_D((n+1)D)) 
\longrightarrow
 \cdot\cdot\cdot.
$$  
So there are infinitely many $n$, such that 
$h^0({\mathcal{O}}_D(nD))>0 $ since 
$h^0(X, {\mathcal{O}}_X(nD))$ grows like $cn^3$
for sufficiently large $n$. 
Let 
$$ N=\{ m\in {\Bbb{N}}, h^0(D, {\mathcal{O}}_D(mD))>0\}.
$$
Let $p$ be the greatest common divisor of $N$, 
then  
$$ h^0(D, {\mathcal{O}}_D(npD)) >0
$$
for all $n\geq 0$.  Thus the line bundle 
${\mathcal{O}}_D(pD)$ determines an effective divisor 
$G$ on $D$ such that 
$${\mathcal{O}}_D(pD)={\mathcal{O}}_D(G).$$  
If ${\mathcal{O}}_D(D)$ is torsion, then by Lecture 3, \cite{I1},
$h^0(D, {\mathcal{O}}_D(nD))\leq 1$ for all $n\geq 0$. 
Since $D$ is big, i.e., 
$h^0(X, {\mathcal{O}}_X(nD))> c n^3$, 
${\mathcal{O}}_D(D)$ is not torsion. 
So we may assume that $D|_D$ determines an effective 
divisor on $D$. We still denote it as $G$. 
By Lemma 2.4, we have 
for all $j\geq 0$ and $i>0$,
$$ \lim_{{\stackrel{\to}{n}}} H^i(D,  \Omega^j_X(nD)|_D)=0.  
$$
In particular, for all $j\geq 0$ and $i>0$,
$$ \lim_{{\stackrel{\to}{n}}} H^i(D,  \Omega^j_X(nD)|_D)=
\lim_{{\stackrel{\to}{n}}} 
H^i(D,  \Omega^j_X\otimes {\mathcal{O}}_D(nD))
$$
$$\quad\quad\quad\quad\quad\quad\quad\quad\quad\quad\quad\quad
=\lim_{{\stackrel{\to}{n}}} 
H^i(D,  \Omega^j_X\otimes {\mathcal{O}}_D(nG))=0.  
$$
Let $S=D-G$, by Lemma 2.1, for all $i>0$ and $j\geq 0$, we have 
$$H^i(S, \Omega^j_X|_S)= 
H^i(S, \Omega^j_D|_S)=0. $$

By the first exact sequence of Lemma 2.3, we have 
 $$ \cdot\cdot\cdot
   \longrightarrow  
 H^1(X,   \lim_{{\stackrel{\to}{n}}}  {\mathcal{O}}_X((n+1)D))
\longrightarrow 
   H^1(D,      \lim_{{\stackrel{\to}{n}}}  
 {\mathcal{O}}_D((n+1)D))
$$
$$
\longrightarrow  
H^2(X, \lim_{{\stackrel{\to}{n}}}  {\mathcal{O}}_X(nD))
\longrightarrow  \cdot\cdot\cdot
$$
Since the direct limit commutes with cohomology
\cite{H1}, Chapter III, Proposition 2.9, and  by Lemma 2.1, 
for all $i>0$, 
$$\lim_{{\stackrel{\to}{n}}}
H^i(X, {\mathcal{O}}_X(nD) )
=0,
$$
we have  
$$ H^i(S, {\mathcal{O}}_S)=
\lim_{{\stackrel{\to}{n}}}
H^i(D, {\mathcal{O}}_D(nD) )=0. 
$$
By \cite{H1}, page 178, Theorem 8.17 or
\cite{GrH}, page 157, ($**$), we have 
$$0\longrightarrow 
{\mathcal{O}}_D(-D)
\longrightarrow 
\Omega^1_X|_D
\longrightarrow 
\Omega^1_D
\longrightarrow 
0.
$$
Tensoring with ${\mathcal{O}}_D((n+1)D)$,
we have
$$0\longrightarrow 
{\mathcal{O}}_D(nD)
\longrightarrow 
\Omega^1_X|_D \otimes {\mathcal{O}}_D((n+1)D)
\longrightarrow 
\Omega^1_D ((n+1)D)
\longrightarrow 
0.
$$
Taking the direct limit, the corresponding long exact sequence gives 
$$H^i(S, \Omega^1_S)=
\lim_{{\stackrel{\to}{n}}}
H^i(D,\Omega^1_D ((n+1)D) )
=0.$$ 
Applying the same procedure to the following short exact sequence
$$ 0\longrightarrow 
\Omega^1_D(-D)
\longrightarrow 
\Omega^2_X|_D
\longrightarrow 
\Omega^2_D
\longrightarrow 
0,
$$ 
we have
$H^i(S, \Omega^2_S)=0$.
So $S$ satisfies the same vanishing condition, i.e., 
for all $i>0$, $j\geq 0$, 
$H^i(S, \Omega^j_S)=0$.

 By Theorem 2.6 in Section 2
 and Lemma 1.8 \cite{Ku},
the boundary divisor 
$G$ on $D$ is connected and $\kappa(G, D)=0$ or $2$.
We know that there are $(-2)$-curves on type (3)
surface in Theorem 2.6. 
 Since $D$
has no exceptional curves, $D$ is not of type (3).
Let $C$ be an irreducible  complete curve on $X$. If $C$ is not contained in
$D$, then $C\cdot D>0$ since $Y$ has no complete curves by Lemma 5 \cite{Zh1}.
If $C$ is contained in $D$ but not in $G$, then  again
$$C\cdot D= C\cdot D|_D=C\cdot G >0$$
since $S$ has no complete curves by Lemma 1.1 \cite{Ku}. 
If $C$ is a component of $G$, then $C\cdot G\geq 0$
since $G$ has no exceptional curves and is connected.
 Therefore  $D$ is nef. 
 We know that 
$S$ is either affine or type (2) surface in Theorem 2.6.
If $S$ is of type (2), then the boundary $D-S$ is an irreducible 
elliptic curve with self-intersection number 0 by Theorem 2.6. This gives us 
$G^2=D^3=0$  which is impossible since $D$ is nef and
$h^0(X, {\mathcal{O}}_X(nD))> c n^3$ 
(Proposition  2.61, \cite{KM}). 
So $S$ must be affine. Therefore 
for any component $C$  of $G$, 
$$ D^3= G^2\geq  C\cdot G>0.
$$
Therefore $D$ is big and nef.

 We will  prove that the linear system $|nD|$
is base  point   free for sufficiently large $n$. 

  Since $S$ is an  affine surface and $G$
contains no  exceptional curves, $G$ is ample, i.e., 
$pD|_D$ is ample on $D$. Thus $D|_D$ is ample on $D$. 
Let $L=D|_D$, then
$$ H^1(D, {\mathcal{O}}_D(nL))=
H^1(D, {\mathcal{O}}_D(nD))=0
$$ 
for sufficiently large $n$.  
Thus  for all  sufficiently large $n$
we have surjective map 
$$ H^1(X, {\mathcal{O}}_X((n-1)D))
\longrightarrow
H^1(X, {\mathcal{O}}_X(nD))
\longrightarrow
0.
$$
Since  
$$ \lim_{{\stackrel{\to}{n}}} 
    H^1(X, {\mathcal{O}}_X(nD))=0,
$$
we  have  $ H^1(X, {\mathcal{O}}_X(nD))=0$ for large $n$.  
Thus we have  an  exact sequence
 $$
 0
\longrightarrow
H^0({X, \mathcal{O}}_X((n-1)D))
\longrightarrow
H^0(X, {\mathcal{O}}_X(nD))
\longrightarrow
H^0(D, {\mathcal{O}}_D(nD))
\longrightarrow
0.
$$
 This implies  that the linear system  $|nD|$ is base  point   free. 
In fact, a point $x\in  X$ is a base point of $|nD|$
if and only if for every section
$s\in H^0({\mathcal{O}}_X(nD))$,  $s(x)=0$  or if and only if  
for every effective divisor 
$E\in |nD|$, $x\in E$.
Suppose that $x$ is a base point of $|nD|$,
then $x\in nD$. Thus $x\in D$.
Since $nD|_D$ is very ample  for large $n$, there is a divisor $F\in |nD|_D|$
such that  $x$ is not a point of $F$. Pull $F$ back 
 to $|nD|$, then there is a divisor $E\in |nD|$, such that  
$x$ is not a point of $E$. This is a contradiction.   
So the linear system $|nD|$
is base point free for sufficiently large $n$. 
\begin{flushright}
 Q.E.D. 
\end{flushright}

{\bf Proof of Theorem 1.3.}
In order to prove Theorem 1.3, we will analyse three  cases: 
$D^3>0$, $D^3< 0$ and $D^3=0$.

\begin{proposition} Let $Y$ be a smooth threefold contained in a 
smooth projective
threefold $X$ such that  $H^0(Y, {\mathcal{O}}_Y )=\Bbb{C}$,
  $H^i(Y, \Omega^j_Y)=0$
for all $j\geq 0$ and $i>0$ and the complement  $D=X-Y$ is a smooth 
projective 
surface. Then   
$D^3$ is not positive. 
\end{proposition}

{\it Proof}. Suppose $D^3>0.$
By the Riemann-Roch formula for surfaces, we have       
$$ h^0({\mathcal{O}}_D(nD))-h^1({\mathcal{O}}_D(nD))+h^2({\mathcal{O}}_D(nD))
=\chi({\mathcal{O}}_D)+\frac{1}{2}n^2D^3-\frac{1}{2}nD^2K_D,
$$  
$$h^0({\mathcal{O}}_D(-nD))-h^1({\mathcal{O}}_D(-nD))+h^2({\mathcal{O}}_D(-nD))
=\chi({\mathcal{O}}_D)+\frac{1}{2}n^2D^3+\frac{1}{2}nD^2K_D.
$$  
So we have 
$$ h^0({\mathcal{O}}_D(nD))+h^2({\mathcal{O}}_D(nD))
\geq \chi({\mathcal{O}}_D)+\frac{1}{2}n^2D^3-\frac{1}{2}nD^2K_D,
$$
$$h^0({\mathcal{O}}_D(-nD))+h^2({\mathcal{O}}_D(-nD))
\geq\chi({\mathcal{O}}_D)+\frac{1}{2}n^2D^3+\frac{1}{2}nD^2K_D.
$$  
Since $D^3>0$ and $(K_D+nD|_D)+(K_D-nD|_D)=2K_D$, either 
$h^0({\mathcal{O}}_D(nD))\longrightarrow \infty$
or 
$h^0({\mathcal{O}}_D(-nD)) \longrightarrow \infty$
as $n\longrightarrow \infty$ and they cannot be big  simultaneously. 
Otherwise, we would have $h^0({\mathcal{O}}_D(2K_D))=\infty$
which is absurd.\\

Case 1. $h^0({\mathcal{O}}_D(nD))\longrightarrow \infty$.

In this case, there is an integer  $n_0>0$ such that  
 for all $n\geq n_0$, $nD|_D$  is an effective divisor on  $D$. 
We may assume that $D|_D=G$ is effective without  
loss of generality. 
From the exact sequence
$$0\longrightarrow 
{\mathcal{O}}_D(nG)
\longrightarrow
{\mathcal{O}}_D((n+1)G) 
\longrightarrow
{\mathcal{O}}_G((n+1)G) 
\longrightarrow
0, 
$$
since $H^2({\mathcal{O}}_G((n+1)G))=0$, we have surjective map 
$$ H^2({\mathcal{O}}_D(nG))\longrightarrow
H^2({\mathcal{O}}_D((n+1)G) )\longrightarrow 0. 
$$
By Lemma 2.4,  
$$  \lim_{{\stackrel{\to}{n}}} 
H^2({\mathcal{O}}_D(nG))
=\lim_{{\stackrel{\to}{n}}} 
H^2({\mathcal{O}}_D(nD))=0.
$$
So for all sufficiently large $n$,
$H^2({\mathcal{O}}_D(nD))=0$.  
By the same argument,  
$H^2({\mathcal{O}}_X(nD))=0$ for $n\gg 0$. 
Since by Lemma 2.5, $H^3({\mathcal{O}}_X(nD))=0$, we have  
$$  \chi({\mathcal{O}}_X(nD))=
h^0({\mathcal{O}}_X(nD))-h^1({\mathcal{O}}_X(nD))
=1-h^1({\mathcal{O}}_X(nD))\leq 1.
$$
On the other hand, since $D^3>0$,  by the Riemann-Rock formula for threefolds 
(\cite{Ful},
page 291), 
we have 
$$ \chi({\mathcal{O}}_X(nD))=
\frac{1}{6}(nD)^3+
\frac{1}{4} c_1\cdot(nD)^2
+\frac{1}{12}(c_1^2+c_2)\cdot nD
+\frac{1}{24}c_1\cdot c_2
\longrightarrow \infty  
$$ 
as $n\longrightarrow \infty .$
Therefore we get a contradiction. 
Thus   $h^0({\mathcal{O}}_D(nD))\longrightarrow \infty$  is not a 
possible case.  \\

Case 2. $h^0({\mathcal{O}}_D(-nD))\longrightarrow \infty$ 
as $n\longrightarrow \infty .$

There is an $n_0>0$ such that for all 
$n\geq n_0$, $-nD|_D$ is effective divisor on 
$D$. Again we may assume that  $G=-D|_D$ is effective. 
Since  $(K_X+D)|_D=K_D$ (Proposition 8.20, Chapter II, \cite{H1}),  
by Lemma 2.4 and Serre duality, we have 
$$ \lim_{{\stackrel{\to}{n}}}  H^2(D, \Omega^3_X((n+1)D)|_D)=
\lim_{{\stackrel{\to}{n}}}  H^2(D, {\mathcal{O}}_X(K_X+(n+1)D)|_D)
$$
$$\quad\quad \quad\quad\quad \quad\quad\quad \quad\quad\quad 
=
\lim_{{\stackrel{\to}{n}}}  H^2(D, {\mathcal{O}}_D(K_D+nD|_D))=0.
$$ 
Take the duality of the   vector  space, we have  
(\cite{H2}, Chapter 3, Section 3)
$$ \lim_{{\stackrel{\gets}{n}}}  
H^0(D, {\mathcal{O}}_D(-nD))
=(\lim_{{\stackrel{\to}{n}}}H^2(D, {\mathcal{O}}_D(K_D+nD)))^*
=0,
$$
where $*$  indicates the  dual  vector space. 
Let $A_n=H^0({\mathcal{O}}_D)$, then $A_n=\Bbb{C}$.
Let $B_n=H^0(D, {\mathcal{O}}_D(-nD))$ and let $C_n$
be their quotient, then we have short exact sequence
$$ 0\longrightarrow 
A_n
\longrightarrow
B_n
\longrightarrow
C_n
\longrightarrow
0. 
$$ 
Since  $-D|_D=G$  is effective, $B_n$ is a subspace of linear 
space  $B_{n'}$ if $n'> n$. The map 
$B_{n'}\rightarrow B_n$ is the natural  restriction map.
So we have the exact sequence of inverse systems 
$$0\longrightarrow 
(A_n)
\longrightarrow
(B_n)
\longrightarrow
(C_n)
\longrightarrow
0. 
$$ 
By \cite{H1}, page 191, 
we have injective map
$$0\longrightarrow  {\Bbb{C}}=
\lim_{{\stackrel{\gets}{n}}}   A_n
\longrightarrow
\lim_{{\stackrel{\gets}{n}}}  
B_n
$$
which is contrary  to 
$$\lim_{{\stackrel{\gets}{n}}}  
B_n=0.
$$
Thus  $h^0({\mathcal{O}}_D(-nD))\longrightarrow \infty$ 
as $n\longrightarrow \infty $ 
is also impossible.

We have seen that $D^3$  is not positive.

\begin{flushright}
 Q.E.D. 
\end{flushright} 

\begin{remark} Since $h^0(X, {\mathcal{O}}_X(nD))=1$,
when $D$  is  not irreducible but nef, we still have $D^3\leq 0$
\cite{KM}, Proposition 2.61.
\end{remark}

\begin{proposition} With the  assumption of Proposition 3.1, 
$D^3$ is not  negative. 
\end{proposition}

{\it Proof}. Suppose $D^3<0.$
If $h^0({\mathcal{O}}_D(nD))=0$ for all $n\gg 0$, 
then we have injective map 
$$ 0\longrightarrow H^1(X, {\mathcal{O}}_X(nD))
\longrightarrow H^1(X, {\mathcal{O}}_X((n+1)D)).
$$
By Lemma 2.4, 
$$ \lim_{{\stackrel{\to}{n}}}H^1(X, {\mathcal{O}}_X(nD))=0.
$$
These two restrictions imply $H^1(X, {\mathcal{O}}_X(nD))=0$
for all $n\gg 0$. Since $H^3(X, {\mathcal{O}}_X(nD))=0$ for all $n\gg 0$,
 the Euler characteristic
$$\chi ({\mathcal{O}}_X(nD))=
1+h^2({\mathcal{O}}_X(nD))\geq 1.
$$
 On the other hand, since  $D^3<0$, by the Riemann-Roch formula 
(\cite{Ful}, page 291), we have
$$  \chi ({\mathcal{O}}_X(nD))= 
 \frac{1}{6} n^3D^3+
\frac{1}{4} c_1\cdot(nD)^2
+\frac{1}{12}(c_1^2+c_2)\cdot nD
+\frac{1}{24}c_1\cdot c_2
$$
$$
\longrightarrow -\infty \quad {\mbox{as}} \quad n\rightarrow \infty.  
$$
This contradicts the fact that $\chi ({\mathcal{O}}_D(nD))$
is positive. Therefore there are infinitely 
many $n$ such that  
$h^0({\mathcal{O}}_D
(nD))>0$. 
 Let $N$, $p$ be the same 
as in the proof of Theorem 1.1,
i.e.,  
$$ N=\{ m\in {\Bbb{N}}, h^0(D, {\mathcal{O}}_D(mD))>0\},
$$
and  $p$ is the greatest common divisor of $N$.
Then   $h^0(D, {\mathcal{O}}_D(pD))>0$. 
Since $D^3\neq 0$,  the line bundle ${\mathcal{O}}_D(pD)$ is not trivial 
and $pD|_D$ defines an effective divisor  on surface $D$.
Thus we may  assume that $D|_D=G$ is an effective divisor on $D$.

Let $S=D-G$, by the same argument as in the proof of Theorem 1.1,
$H^i(S, \Omega^j_S)=0$ for all $i>0$ and $j\geq 0$
since  
$$\lim_{{\stackrel{\to}{n}}}
H^i(D,  \Omega^j_X(nD)|_D)=0.
$$  
By Theorem 2.6 in Section 2 and Lemma 1.8  \cite{Ku}, $\kappa({G, D})=0$ or $2$. 
If $\kappa({G, D})=2$, then by Theorem 2.7, for all $n\gg 0$,
$$  h^0(D, {\mathcal{O}}_D(nG))=a_2n^2+a_1n+a_0+\lambda(n),
$$
where $\lambda(n)$ is periodic. By Theorem 2.8, 
since $\kappa({G, D})=2$, 
there is a constant $c>0$ such that 
$ h^0(D, {\mathcal{O}}_D(nG))> cn^2$. Thus $a_2>0$. 
By Lemma 2.1, 
$$ \lim_{{\stackrel{\to}{n}}}H^0(D, {\mathcal{O}}_D(nG))
=\lim_{{\stackrel{\to}{n}}}H^0(D, {\mathcal{O}}_D(nD))
=H^0(S, {\mathcal{O}}_S) \neq 0. 
$$
In fact, since $S$ is affine, $h^0(S, {\mathcal{O}}_S)=\infty$.
By Lemma  2.3, we have 
 $$ 0\longrightarrow \lim_{{\stackrel{\to}{n}}}
H^0(X, {\mathcal{O}}_X(nD-D))
\longrightarrow \lim_{{\stackrel{\to}{n}}}
H^0(X, {\mathcal{O}}_X(nD))
$$
$$
\longrightarrow \lim_{{\stackrel{\to}{n}}}
H^0(X, {\mathcal{O}}_D(nD))
\longrightarrow 
0.
$$
By Lemma 2.1, 
$$ \lim_{{\stackrel{\to}{n}}}
H^0(X, {\mathcal{O}}_X(nD))=H^0(Y, {\mathcal{O}}_Y)
=\Bbb{C}
$$
and
$$\lim_{{\stackrel{\to}{n}}}
H^0(X, {\mathcal{O}}_X(nD-D))=H^0(Y, {\mathcal{O}}_X(-D)|_Y)
=\Bbb{C}. 
$$
Thus the above exact sequence is
$$ 0\longrightarrow 
{\Bbb{C}} \longrightarrow 
{\Bbb{C}}
\longrightarrow 
H^0(S, {\mathcal{O}}_S)
\longrightarrow 
0.
$$
This is impossible since $h^0(S, {\mathcal{O}}_S)=\infty$. 
So $\kappa(G, D)\neq 2$.

If $\kappa(G, D)=0$, since
   ${\mathcal{O}}_D(D)\neq  {\mathcal{O}}_D $ defines an effective 
divisor  on 
$D$, by Lemma 2.1, 2.4,
and  Lemma 1.8\cite{Ku} we have 
 $$ \lim_{{\stackrel{\to}{n}}}
H^0(D, {\mathcal{O}}_D(nD))=
H^0(S, {\mathcal{O}}_S)
=\Bbb{C}.
$$
Again  by Lemma 2.1 and 2.3, we have 
$$
0\longrightarrow \lim_{{\stackrel{\to}{n}}}
H^0(X, {\mathcal{O}}_X(nD-D))
\longrightarrow \lim_{{\stackrel{\to}{n}}}
H^0(X, {\mathcal{O}}_X(nD))
$$
$$
\longrightarrow \lim_{{\stackrel{\to}{n}}}
H^0(D, {\mathcal{O}}_D(nD))
\longrightarrow 
0,
$$
which is impossible since three direct limits 
are $\Bbb{C}$. Therefore $D^3$  is not negative.

\begin{flushright}
 Q.E.D. 
\end{flushright}

\begin{proposition} Under the condition of Proposition 3.1, 
 $D^3=0$  and 
$D|_D$ is numerically equivalent to 0. 
\end{proposition}

{\it Proof.}
Let $H$ be an ample divisor on $X$. Then when restricted on $D$,
$L=H|_D$ is still ample by Proposition 4.1 \cite{H2}. We may 
assume that $L$ is smooth
by Bertini's Theorem (Chapter 2, Theorem 8.8 \cite{H1}
or  Theorem 4.21, \cite{ Uen1}). 
Let $G$ be the divisor determined by the line bundle 
${\mathcal{O}}_D(D)$.

If $L\cdot G \neq 0$,  then 
there exists an $n_0>0$ such that either
the linear system  $|n_0G|$ is nonempty
or $|-n_0G|$ is nonempty. If  $|n_0G|$ is nonempty,
we may assume that $G$ is effective.   
Let $S=D-G$. 
 By Lemma 2.3,  we have 
$$\lim_{{\stackrel{\to}{n}}}
H^0(X, {\mathcal{O}}_D(nD))=H^0(S, {\mathcal{O}}_S).
$$
By  Theorem 2.8, 
$$ h^0(S, {\mathcal{O}}_S)=
{\mbox{dim}}_{\Bbb{C}}
\lim_{{\stackrel{\to}{n}}}
H^0(X, {\mathcal{O}}_D(nD))\neq 0.
$$
By  the first exact sequence in Lemma 2.3, 
we have exact sequence 
$$0\longrightarrow  
{\Bbb{C}}
\longrightarrow  
{\Bbb{C}}
\longrightarrow  
H^0(S, {\mathcal{O}}_S)
\longrightarrow  
0.
$$ 
This is impossible. 

   If $|-n_0G|$ is nonempty, then 
we may assume that 
 $-G$ is an  effective divisor on $D$. By
 the same argument as in case 2, Proposition 3.1,  we have 
 $$ \lim_{{\stackrel{\gets}{n}}}  
H^0(D, {\mathcal{O}}_D(-nD))
=(\lim_{{\stackrel{\to}{n}}}H^2(D, {\mathcal{O}}_D(K_D+nD)))^*
=0,
$$
which is again not  possible
since 
  $$ {\Bbb{C}}=\lim_{{\stackrel{\gets}{n}}} H^0(D, {\mathcal{O}}_D)
\hookrightarrow \lim_{{\stackrel{\gets}{n}}}  
H^0(D, {\mathcal{O}}_D(-nD)).
$$

We have seen that the only  possible case is 
$L\cdot G=0$. Since $G^2=D^3=0$, by the 
Hodge Index Theorem,
$G$ is numerically  equivalent  to 
0. 
If  $G=0$, then 
$$(\lim_{{\stackrel{\to}{n}}}H^2(D, {\mathcal{O}}_D(K_D+nD)))^*
=\lim_{{\stackrel{\gets}{n}}}  
H^0(D, {\mathcal{O}}_D(-nD)) =\Bbb{C}
$$
which contradicts Lemma 2.4.
So we have $G\equiv 0$ but $G\neq 0$.

\begin{flushright}
 Q.E.D. 
\end{flushright}

\begin{proposition} Under the condition of Proposition 3.1, 
If ${\mathcal{O}}_D(D)$ is not torsion,  then $q=h^1(X, {\mathcal{O}}_X)=0$,
$\frac{1}{2}(c_1^2+c_2)\cdot D=\chi({\mathcal{O}}_D)=0$,
 $\chi ({\mathcal{O}}_X)>0$ and $K_X$ is not nef.
\end{proposition}  
{\it Proof.}
If there is  an $n_0$ such that $h^0(D, {\mathcal{O}}_D(n_0D))\neq 0$
(or  $h^0(D, {\mathcal{O}}_D(-n_0D))\neq 0$), then $\kappa(D, D|_D)\geq 0$.
By \cite{I1}, page 34-35, 
 there are infinitely many 
$n$ such that $h^0(D, {\mathcal{O}}_D(nD))\neq 0$
(or  $h^0(D, {\mathcal{O}}_D(-nD))\neq 0$).
Then  $D|_D$  (or $-D|_D$)  defines an effective divisor
on $D$ since ${\mathcal{O}}_D(D)$ is not torsion. By  Lemma 2.3,
we can get a
contradiction
$$ 0 \longrightarrow {\Bbb{C}}
\longrightarrow {\Bbb{C}}
\longrightarrow H^0(S, {\mathcal{O}}_S)
\longrightarrow  0.
$$
 So 
$h^0(D, {\mathcal{O}}_D(nD))=h^0(D, {\mathcal{O}}_D(-nD))=0$
for all $n>0$.

Since $h^0(D, {\mathcal{O}}_D(nD))=0$ for all $n>0$, from
the exact sequence
$$0\longrightarrow 
{\mathcal{O}}_X((n-1)D)
\longrightarrow
{\mathcal{O}}_X(nD) 
\longrightarrow
{\mathcal{O}}_D(nD) 
\longrightarrow
0, 
$$
we have injective map from $H^1(X, {\mathcal{O}}_X(nD))$
to $H^1(X, {\mathcal{O}}_X((n+1)D))$. Since the direct limit is 0,
$H^1(X, {\mathcal{O}}_X(nD))=0$ for sufficiently large $n$.
By the injectivity, $H^1(X, {\mathcal{O}}_X((n-1)D))=0$ for
all $n>0$, i.e., 
$H^1(X, {\mathcal{O}}_X(nD))=H^1(X, {\mathcal{O}}_X)=0$. 
Thus for all $n > 0$, by Lemma 2.5,  we have exact sequence
$$0\longrightarrow  H^1(D, {\mathcal{O}}_D((n+1)D))
\longrightarrow  H^2(X, {\mathcal{O}}_X(nD))
\longrightarrow  H^2(X, {\mathcal{O}}_X((n+1)D))
$$
$$
\longrightarrow  H^2(D, {\mathcal{O}}_D((n+1)D))
 \longrightarrow
0.
$$
This gives us 
$$  h^2( {\mathcal{O}}_X((n+1)D))-
h^2( {\mathcal{O}}_X(nD))=
h^2( {\mathcal{O}}_D((n+1)D))
-h^1( {\mathcal{O}}_D((n+1)D)).
$$
By Riemann-Roch formulas for surfaces and
threefolds, we have  
$$  \frac{1}{12} (c_1^2+c_2)\cdot D = 
\chi ({\mathcal{O}}_D). 
$$
Since  for all  $n\gg 0$, $h^1({\mathcal{O}}_X(nD))=h^3({\mathcal{O}}_X(nD))=0$
and $D^3=D^2\cdot K_D=0$,
we have 
$\chi({\mathcal{O}}_X(nD))=1+ h^2( {\mathcal{O}}_X(nD))\geq 1$.
So
$\frac{1}{12} (c_1^2+c_2)\cdot D = \chi ({\mathcal{O}}_D)\geq 0. 
$

If  $n\gg 0$, $h^1({\mathcal{O}}_D(nD))=0$, then 
by the above exact sequence, we have  
$$0
\longrightarrow  H^2(X, {\mathcal{O}}_X(nD))
\longrightarrow  H^2(X, {\mathcal{O}}_X((n+1)D))
\longrightarrow  H^2(D, {\mathcal{O}}_D((n+1)D))
 \longrightarrow
0,
$$
It is easy to see that for all $n\geq 0$, $h^2(X, {\mathcal{O}}_X(nD))=0$
and $h^2( {\mathcal{O}}_D((n+1)D))=h^1( {\mathcal{O}}_D((n+1)D))=0$. 
Thus  $\chi ({\mathcal{O}}_D)=0$ and the proposition follows.

So the remaining case is:  for infinitely many $n>0$, 
$h^1({\mathcal{O}}_D(nD))>0$. 
We claim that $\chi ({\mathcal{O}}_D)=0$. If not, 
then $\chi ({\mathcal{O}}_D)>0$.
Since for all $n>0$, $h^0(D, {\mathcal{O}}_D(-nD))=
h^2(D, {\mathcal{O}}_D(K_D+nD))=0$, we have surjective map from 
$H^2(X, {\mathcal{O}}_X(K_X+(n-1)D))$
to $H^2(X, {\mathcal{O}}_X(K_X+nD))$. By Lemma 2.4, 
$$\lim_{{\stackrel{\to}{n}}}H^2(X, {\mathcal{O}}_X(K_X+nD))=0.
$$
So for $n\gg 0$, $H^2(X, {\mathcal{O}}_X(K_X+nD))=0$.

Applying  Riemann-Roch to the sheaf ${\mathcal{O}}_X(K_X+nD)$, we have 
$$ h^0(X, {\mathcal{O}}_X(K_X+nD))\geq 
(\frac{1}{12} (c_1^2+c_2)\cdot D )n>0.
$$
Since $h^0(X, {\mathcal{O}}_X(nD))=1$ for all $n\geq 0$, the above 
inequality is not true. To see this, let $D_m$ be an effective divisor
in the linear system $|K_X+ mD|$ and 
$D_{m+n}$ be an effective divisor in the linear system 
$|K_X+(m+n)D|$, then  $D_m+nD$ is linearly equivalent to $D_{m+n}$.
Let $f$ be the nonconstant rational function on $X$ such that
div$f+D_{m+n}=D_m+nD$, then  div$f+D_{m+n}-D_m=nD>0$. This gives 
$f\in H^0(X, {\mathcal{O}}_X(D_{m+n}-D_m))$ 
and $D_{m+n}-D_m$ is linearly equivalent to $nD$. 
Since $h^0(X, {\mathcal{O}}_X(K_X+(m+n)D))\geq c(m+n)$, there are 
at least 2  linearly independent   $f$ in 
$H^0(X, {\mathcal{O}}_X(D_{m+n}-D_m))$. 
  It   contradicts to $H^0(X, {\mathcal{O}}_X(D_{m+n}-D_m))=
H^0(X, {\mathcal{O}}_X(nD)=\Bbb{C}$ for all $n\geq 0$. 
Thus
$h^0(X, {\mathcal{O}}_X(K_X+nD))$ cannot be greater than $cn$ for $c>0$.
So  
$$  \frac{1}{12} (c_1^2+c_2)\cdot D = 
\chi ({\mathcal{O}}_D)=0. 
$$
By Riemann-Roch,  $h^2(X, {\mathcal{O}}_X(nD))=h^2(X, {\mathcal{O}}_X)-
h^3(X, {\mathcal{O}}_X)$ is a constant 
and $h^2( {\mathcal{O}}_D(nD))
=h^1( {\mathcal{O}}_D(nD))$  for all $n\geq 0$. 
Since $h^2(X, {\mathcal{O}}_X)-
h^3(X, {\mathcal{O}}_X)\geq 0$ and $q=0$, 
$\chi ({\mathcal{O}}_X)\geq 1$. By a theorem of Miyaoka
\cite {Mi}, $K_X$ is not nef.

\begin{flushright}
 Q.E.D. 
\end{flushright}

\begin{proposition} Under the condition of Proposition 3.1, 
if ${\mathcal{O}}_D(D)$ is torsion, then $h^1(X, {\mathcal{O}}_X(nD))$
are bounded for all $n$ and 
$ \frac{1}{12} (c_1^2+c_2)\cdot D = 
\chi ({\mathcal{O}}_D)\geq 0$. 
\end{proposition}
{\it Proof. }  Let  $p$ be the least natural number such that 
${\mathcal{O}}_D(pD)={\mathcal{O}}_D$. 
Let $G$ be the divisor defined by $D|_D$, then   
$\kappa (G, D)=0$. 
By Iitaka \cite{I1}, page 35, 
$h^0(D, {\mathcal{O}}_D(nD))=0$ if $p$ does not divide $n$. 
Lemma 2.1-2.5  in \cite{Ku} still hold in our case, 
so  $h^1(X, {\mathcal{O}}_X(nD))$ are bounded for all $n\geq 0$.
Since $h^3(X, {\mathcal{O}}_X(nD))=0$ and $D^3=D^2\cdot K_X=0$,
by Riemann-Roch, we have 
$$\chi ({\mathcal{O}}_X(nD))=1-h^1(X, {\mathcal{O}}_X(nD))
+h^2(X, {\mathcal{O}}_X(nD))=\frac{1}{12} (c_1^2+c_2) nD
+\frac{1}{24}c_1c_2.$$ 
So  $\frac{1}{12} (c_1^2+c_2) D\geq 0$ since 
$h^1(X, {\mathcal{O}}_X(nD))$ are bounded for all $n\geq 0$.

 By the short exact sequence
$$ 0
\longrightarrow   {\mathcal{O}}_X((np-1)D)
\longrightarrow   {\mathcal{O}}_X(npD)
\longrightarrow  {\mathcal{O}}_D(npD)={\mathcal{O}}_D
 \longrightarrow
0,  
$$
we have  
$$\chi ({\mathcal{O}}_D)=
\chi ({\mathcal{O}}_X(npD))-\chi ({\mathcal{O}}_X((np-1)D))
=\frac{1}{12} (c_1^2+c_2) D\geq 0.
$$

\begin{flushright}
 Q.E.D. 
\end{flushright}

  We  have proved  Theorem 1.3.

By the standard classification theory of surfaces, 
if ${\mathcal{O}}_D(D)$ is not torsion,  $D$ might be
one of the following surfaces: (1) ruled surfaces over an elliptic curve;
(2) bi-elliptic surfaces; (3) Abelian surfaces; (4) Elliptic surfaces.

\end{document}